\documentclass{amsart}

\usepackage{palatino}	

\theoremstyle{plain}
\newtheorem{Theorem}{Theorem}
\newtheorem{Lemma}[Theorem]{Lemma}

\newtheorem{Proposition}[Theorem]{Proposition}
\newtheorem{Problem}[Theorem]{Problem}

\theoremstyle{remark}

\newtheorem{Remarks}[Theorem]{Remarks}

\theoremstyle{definition}
\newtheorem{Definition}[Theorem]{Definition}
\newtheorem{Notation}[Theorem]{Notation}
\newtheorem{Example}[Theorem]{Example}
\newtheorem{Proof}{Proof}

\usepackage[black-square]{QED}
\usepackage{amssymb}
\usepackage{mathrsfs} 
\usepackage{amsmath} 
 
\def\diatop[#1|#2]{{\setbox1=\hbox{{#1{}}}\setbox2=\hbox{{#2{}}}%
                    \dimen0=\ifdim\wd1>\wd2\wd1\else\wd2\fi%
                    \dimen1=\ht2\advance\dimen1by-1ex%
                    \setbox1=\hbox to1\dimen0{\hss#1\hss}%
                    \rlap{\raise1\dimen1\box1}%
                    \hbox to1\dimen0{\hss#2\hss}}}%
 
\font\ipatwelverm=wsuipa12
\def\ipa{\ipatwelverm}
 
\def\inva{\edef\next{\the\font}\ipa\char'000\next}%
\def\scripta{\edef\next{\the\font}\ipa\char'001\next}%
\def\nialpha{\edef\next{\the\font}\ipa\char'002\next}%
\def\invscripta{\edef\next{\the\font}\ipa\char'003\next}%
\def\invv{\edef\next{\the\font}\ipa\char'004\next}%
 
\def\crossb{\edef\next{\the\font}\ipa\char'005\next}%
\def\barb{\edef\next{\the\font}\ipa\char'006\next}%
\def\slashb{\edef\next{\the\font}\ipa\char'007\next}%
\def\hookb{\edef\next{\the\font}\ipa\char'010\next}%
\def\nibeta{\edef\next{\the\font}\ipa\char'011\next}%
 
\def\slashc{\edef\next{\the\font}\ipa\char'012\next}%
\def\curlyc{\edef\next{\the\font}\ipa\char'013\next}%
\def\clickc{\edef\next{\the\font}\ipa\char'014\next}%
 
\def\crossd{\edef\next{\the\font}\ipa\char'015\next}%
\def\bard{\edef\next{\the\font}\ipa\char'016\next}%
\def\slashd{\edef\next{\the\font}\ipa\char'017\next}%
\def\hookd{\edef\next{\the\font}\ipa\char'020\next}%
\def\taild{\edef\next{\the\font}\ipa\char'021\next}%
\def\dz{\edef\next{\the\font}\ipa\char'022\next}%
\def\eth{\edef\next{\the\font}\ipa\char'023\next}%
\def\scd{\edef\next{\the\font}\ipa\char'024\next}%
 
\def\schwa{\edef\next{\the\font}\ipa\char'025\next}%
\def\er{\edef\next{\the\font}\ipa\char'026\next}%
\def\reve{\edef\next{\the\font}\ipa\char'027\next}%
\def\niepsilon{\edef\next{\the\font}\ipa\char'030\next}%
\def\revepsilon{\edef\next{\the\font}\ipa\char'031\next}%
\def\hookrevepsilon{\edef\next{\the\font}\ipa\char'032\next}%
\def\closedrevepsilon{\edef\next{\the\font}\ipa\char'033\next}%
 
\def\scriptg{\edef\next{\the\font}\ipa\char'034\next}%
\def\hookg{\edef\next{\the\font}\ipa\char'035\next}%
\def\scg{\edef\next{\the\font}\ipa\char'036\next}%
\def\nigamma{\edef\next{\the\font}\ipa\char'037\next}
\def\ipagamma{\edef\next{\the\font}\ipa\char'040\next}%
\def\babygamma{\edef\next{\the\font}\ipa\char'041\next}%
 
\def\hv{\edef\next{\the\font}\ipa\char'042\next}%
\def\crossh{\edef\next{\the\font}\ipa\char'043\next}%
\def\hookh{\edef\next{\the\font}\ipa\char'044\next}%
\def\hookheng{\edef\next{\the\font}\ipa\char'045\next}%
\def\invh{\edef\next{\the\font}\ipa\char'046\next}%
 
\def\bari{\edef\next{\the\font}\ipa\char'047\next}%
\def\dlbari{\edef\next{\the\font}\ipa\char'050\next}
\def\niiota{\edef\next{\the\font}\ipa\char'051\next}%
\def\sci{\edef\next{\the\font}\ipa\char'052\next}%
\def\barsci{\edef\next{\the\font}\ipa\char'053\next}
 
\def\invf{\edef\next{\the\font}\ipa\char'054\next}%
 
\def\tildel{\edef\next{\the\font}\ipa\char'055\next}%
\def\barl{\edef\next{\the\font}\ipa\char'056\next}%
\def\latfric{\edef\next{\the\font}\ipa\char'057\next}%
\def\taill{\edef\next{\the\font}\ipa\char'060\next}%
\def\lz{\edef\next{\the\font}\ipa\char'061\next}%
\def\nilambda{\edef\next{\the\font}\ipa\char'062\next}%
\def\crossnilambda{\edef\next{\the\font}\ipa\char'063\next}%
 
\def\labdentalnas{\edef\next{\the\font}\ipa\char'064\next}%
\def\invm{\edef\next{\the\font}\ipa\char'065\next}%
\def\legm{\edef\next{\the\font}\ipa\char'066\next}%
 
\def\nj{\edef\next{\the\font}\ipa\char'067\next}%
\def\eng{\edef\next{\the\font}\ipa\char'070\next}%
\def\tailn{\edef\next{\the\font}\ipa\char'071\next}%
\def\scn{\edef\next{\the\font}\ipa\char'072\next}%
 
\def\clickb{\edef\next{\the\font}\ipa\char'073\next}%
\def\baro{\edef\next{\the\font}\ipa\char'074\next}%
\def\openo{\edef\next{\the\font}\ipa\char'075\next}%
\def\niomega{\edef\next{\the\font}\ipa\char'076\next}%
\def\closedniomega{\edef\next{\the\font}\ipa\char'077\next}%
\def\oo{\edef\next{\the\font}\ipa\char'100\next}%
 
\def\barp{\edef\next{\the\font}\ipa\char'101\next}%
\def\thorn{\edef\next{\the\font}\ipa\char'102\next}%
\def\niphi{\edef\next{\the\font}\ipa\char'103\next}%
 
\def\flapr{\edef\next{\the\font}\ipa\char'104\next}%
\def\legr{\edef\next{\the\font}\ipa\char'105\next}%
\def\tailr{\edef\next{\the\font}\ipa\char'106\next}%
\def\invr{\edef\next{\the\font}\ipa\char'107\next}%
\def\tailinvr{\edef\next{\the\font}\ipa\char'110\next}%
\def\invlegr{\edef\next{\the\font}\ipa\char'111\next}%
\def\scr{\edef\next{\the\font}\ipa\char'112\next}%
\def\invscr{\edef\next{\the\font}\ipa\char'113\next}%
 
\def\tails{\edef\next{\the\font}\ipa\char'114\next}%
\def\esh{\edef\next{\the\font}\ipa\char'115\next}%
\def\curlyesh{\edef\next{\the\font}\ipa\char'116\next}%
\def\nisigma{\edef\next{\the\font}\ipa\char'117\next}%
 
\def\tailt{\edef\next{\the\font}\ipa\char'120\next}%
\def\tesh{\edef\next{\the\font}\ipa\char'121\next}%
\def\clickt{\edef\next{\the\font}\ipa\char'122\next}%
\def\nitheta{\edef\next{\the\font}\ipa\char'123\next}%
 
\def\baru{\edef\next{\the\font}\ipa\char'124\next}%
\def\slashu{\edef\next{\the\font}\ipa\char'125\next}%
\def\niupsilon{\edef\next{\the\font}\ipa\char'126\next}%
\def\scu{\edef\next{\the\font}\ipa\char'127\next}%
\def\barscu{\edef\next{\the\font}\ipa\char'130\next}%
 
\def\scriptv{\edef\next{\the\font}\ipa\char'131\next}%
 
\def\invw{\edef\next{\the\font}\ipa\char'132\next}%
 
\def\nichi{\edef\next{\the\font}\ipa\char'133\next}%
 
\def\invy{\edef\next{\the\font}\ipa\char'134\next}%
\def\scy{\edef\next{\the\font}\ipa\char'135\next}%
 
\def\curlyz{\edef\next{\the\font}\ipa\char'136\next}%
\def\tailz{\edef\next{\the\font}\ipa\char'137\next}%
\def\yogh{\edef\next{\the\font}\ipa\char'140\next}%
\def\curlyyogh{\edef\next{\the\font}\ipa\char'141\next}%
 
\def\glotstop{\edef\next{\the\font}\ipa\char'142\next}%
\def\revglotstop{\edef\next{\the\font}\ipa\char'143\next}%
\def\invglotstop{\edef\next{\the\font}\ipa\char'144\next}%
\def\ejective{\edef\next{\the\font}\ipa\char'145\next}%
\def\reveject{\edef\next{\the\font}\ipa\char'146\next}%
 
\def\upt{\edef\next{\the\font}\ipa\char'154\next}
\def\downt{\edef\next{\the\font}\ipa\char'155\next}%
\def\leftt{\edef\next{\the\font}\ipa\char'156\next}%
\def\rightt{\edef\next{\the\font}\ipa\char'157\next}%
 
\def\upp{\edef\next{\the\font}\ipa\char'164\next}
\def\downp{\edef\next{\the\font}\ipa\char'165\next}%
\def\leftp{\edef\next{\the\font}\ipa\char'166\next}%
\def\rightp{\edef\next{\the\font}\ipa\char'167\next}%
 
\def\stress{\edef\next{\the\font}\ipa\char'150\next}
\def\secstress{\edef\next{\the\font}\ipa\char'151\next}
 
\def\syllabic{\edef\next{\the\font}\ipa\char'152\next}
 
\def\corner{\edef\next{\the\font}\ipa\char'153\next}%
 
\def\halflength{\edef\next{\the\font}\ipa\char'160\next}
\def\length{\edef\next{\the\font}\ipa\char'161\next}
 
\def\underdots{\edef\next{\the\font}\ipa\char'162\next}%
 
\def\ain{\edef\next{\the\font}\ipa\char'163\next}
 
\def\overring{\edef\next{\the\font}\ipa\char'170\next}%
\def\underring{\edef\next{\the\font}\ipa\char'171\next}%
 
\def\open{\edef\next{\the\font}\ipa\char'172\next}%
 
\def\midtilde{\edef\next{\the\font}\ipa\char'173\next}%
\def\undertilde{\edef\next{\the\font}\ipa\char'174\next}%
 
\def\underwedge{\edef\next{\the\font}\ipa\char'175\next}%
 
\def\polishhook{\edef\next{\the\font}\ipa\char'176\next}%

\usepackage{rotating} 

\usepackage[colorlinks=true, linkcolor=blue, pdfstartview=FitH, pdfauthor={Lukasz Grabarek}]{hyperref}

\newcommand{\Nat}{\mathbb{N}}

\newcommand{\Rea}{\mathbb{R}}
\newcommand{\Com}{\mathbb{C}}

\newcommand{\R}{\mbox{Re \!}}

\newcommand{\sm}{\smallsetminus}

\newcommand{\LP}{$\mathscr{L}$-$\mathscr{P}$} 
\newcommand{\LPP}{$\mathscr{L}$-$\mathscr{P}^+$ } 
\newcommand{\okinaA}{\scalebox{.7}{\raisebox{.45 cm}{\rotatebox{180}{\reflectbox{\ain}}}}}

\begin{document}
\title{A New Class of Non-Linear Stability Preserving Operators}
\author{Lukasz Grabarek}
\address{Department of Mathematics,
University of Hawai\okinaA i at M\=anoa, Honolulu, HI 96822}
\email{lukasz@math.hawaii.edu}

\date{\today}

\begin{abstract}

We extend Br\"and\'en's recent proof of a conjecture of Stanley and describe a new class of non-linear operators that preserve weak Hurwitz stability and the Laguerre-P\'olya class.
\end{abstract}


\maketitle

\section{Introduction}
A real entire function $\psi(x)$ is said to belong to the Laguerre-P\'olya class, denoted \LP, if $\psi(x)$ is the uniform limit, on compact subsets of $\Com$, of real polynomials all of whose zeros are real.  $\psi(x) \in$ \LP \; if and only if $\psi(x)$ can be expressed in the form
\begin{equation}
\label{LP}\large{
\psi(x) = cx^me^{-\alpha x^2+\beta x} \prod_{k=1}^{\infty}\left(1  + \frac{x}{x_k}\right)e^{-\frac{x}{x_k}} } \;,
\end{equation}
where $m \in \Nat, \, \alpha, \beta, c, x_k \in \Rea, \, \alpha \geq 0$, and $\sum_{j=1}^{\infty} \frac{1}{x_k^2} < \infty$.  An important subclass of the Laguerre-P\'olya class, denoted \LPP, consists of precisely those $\varphi(x) \in$ \LP \; whose Taylor coefficients are non-negative.  $\varphi(x) \in$ \LPP if and only if $\varphi(x)$ can be expressed in the form
\begin{equation}
\label{LP+} \large{
\varphi(x) = cx^me^{\sigma x} \prod_{k=1}^{\infty}\left(1  + \frac{x}{x_k}\right) } \;,
\end{equation}
where $m \in \Nat, \, \sigma, c, x_k \in \Rea, \, \sigma, c \geq 0, \, x_k >0$, and $\sum_{j=1}^{\infty} \frac{1}{x_k} < \infty$.  For various properties of the Laguerre-P\'olya class refer to \cite[Ch. VIII]{levin}, \cite{O63}, \cite{PS}, and the references contained therein.  In particular, for $\psi(x) = \sum_{k=0}^{\infty} \frac{\gamma_k}{k!} x^k\in$ \LP \,: 
\begin{enumerate}
\item If $\gamma_k \geq 0$ for $k = 0, 1, 2, \ldots$ \! , then $\psi(x) \in$ \LPP \! ;
\item The \textit{Tur\'an inequalities} hold: $\gamma_k^2 - \gamma_{k-1}\gamma_{k+1} \geq 0$ for $k = 1, 2, \ldots$ \! ;
\item The \textit{Laguerre inequalities} hold: $(\psi^{(p)}(x))^2 - \psi^{(p-1)}(x)\psi^{(p+1)}(x)\geq 0$ for each $p = 1, 2, \ldots$ , and all $x \in \Rea$ \! ;
\item $\psi^{(k)}(x) \in$ \LP \; for $k = 1, 2, \ldots$ .
\end{enumerate}

\vspace{0.5 cm}
The Laguerre inequalities generalize to a system of inequalities which characterize the Laguerre-P\'olya class.  

\begin{Theorem}[\cite{CC89}, \cite{CC02}, \cite{Patrick}]
\label{TLaguerre}
Let $\varphi(x) = \sum_{k=0}^{\infty} \frac{\gamma_k}{k!}x^k$ be a real entire function and define
\begin{equation}\label{Laguerre}
L_p(\varphi(x)) := \sum_{j=0}^{2p} \frac{(-1)^{p+j}}{(2p)!} {2p \choose j} \varphi^{(j)}(x)\varphi^{(2p-j)}(x)\;, \end{equation} 
where $x \in \Rea$ and $p = 0, 1, 2, \ldots$ . Then $\varphi(x) \in$ \LP \; if and only if for all $x \in \Rea$ and all $p = 0, 1, 2, \ldots$
\[ L_p(\varphi(x)) \geq 0 \; .\]
\end{Theorem}

In their 1989 study of the relationship between the Laguerre and the Tur\'an inequalities, Craven and Csordas have posed  the following problem:
\begin{Problem}[\cite{CC89}]\label{CS} Classify the functions 

\begin{equation} \psi(x) = \sum_{k=1}^{\omega} \frac{\gamma_k}{k!} x^k  \in \mbox{\LP} \;, \end{equation}
where $\gamma_k \geq 0$ and $0\leq \omega \leq \infty$, for which the functions 

\begin{equation} f(x) := \sum_{k=0}^{\infty} \frac{\gamma_{k}^2 - \gamma_{k-1}\gamma_{k+1}}{k!} x^k \in \mbox{\LP} \;. \end{equation}
\end{Problem}

\vspace{0.5 cm}
This problem naturally leads to the following result, conjectured independently by Stanley, McNamara-Sagan, and Fisk and proved by Br\"and\'en (\cite{BR10}).
\begin{Theorem}[\cite{BR10}]\label{Stanley} If the zeros of the real polynomial $\psi(x) = \sum_{k=0}^{n} a_kz^k$ are all real and negative, then the zeros of the polynomial
\begin{equation}
\label{Stanley1}
\sum_{k=0}^n(a_k^2 - a_{k-1}a_{k+1})z^k \, , \; \mbox{where } a_{-1} := 0 \mbox{ and } a_{n+1} := 0 \;, \end{equation}
are all real and negative.
\end{Theorem}

\vspace{0.5 cm}
The coefficients in (\ref{Stanley1}) are obtained by means of the non-linear operator $a_k \mapsto a_k^2 - a_{k-1}a_{k+1}$.  We extend Theorem \ref{Stanley} to a class of non-linear operators that take real polynomials with all real negative zeros into polynomials of the same type, and include the non-linear operator $a_k \mapsto a_k^2 - a_{k-1}a_{k+1}$.

\section{Preliminary Results}
\subsection*{A Class of Non-Linear Operators}

\begin{Notation}\label{not1} We will always explicitly write $a_k = \frac{\gamma_k}{k!}$ to avoid confusion (cf. \S 4) and distinguish between the coefficients $a_k$ of $\psi(x)= \sum_{k=0}^{\infty} a_k x^k$ and the coefficients $\frac{\gamma_k}{k!}$ of $\varphi(x)= \sum_{k=0}^{\infty} \frac{\gamma_k}{k!} x^k$.  We will follow the convention on the integer index $k$ (cf. Theorem \ref{Stanley} and Problem \ref{CS}) that $a_k = 0$, whenever $k$ is negative for transcendental entire functions $\psi(x)$, and if $\psi(x)$ is a polynomial, then we will set $a_k = 0$, whenever $k \not\in \{0, 1, 2, \ldots, \deg \psi(x)\}$.
\end{Notation}

\begin{Lemma}\label{Lk} 
Let $\varphi(x) = \sum_{k=0}^{\infty} \frac{\gamma_k}{k!}x^k$ be a real entire function and, for a positive integer $p$, let $L_p\left(\varphi^{(k)}(x) \right)$ be defined as in Theorem \ref{TLaguerre}.  Then, 
\begin{equation}
\left. \frac{(2p)!}{2} \cdot L_p\left( \varphi^{(k)}(x) \right) \right|_{x=0} = {2p-1 \choose p} \gamma_{k+p}^2 + \sum_{j=1}^p (-1)^j {2p \choose p-j} \gamma_{k+p-j}\gamma_{k+p+j} \; .
\end{equation}
\end{Lemma}

\begin{Proof}
Rewriting equation (\ref{Laguerre}) with $\varphi^{(k)}(x)$, yields
\begin{equation}\label{Laguerre2}
L_p\left( \varphi^{(k)}(x) \right) = \sum_{j=0}^{2p} \frac{(-1)^{p+j}}{(2p)!} {2p \choose j} \varphi^{(k+j)}(x)\varphi^{(2p+k-j)}(x) \; .
\end{equation}
For a fixed positive integer $p$, the coefficient of $\gamma^2_{k+p}$ is obtained by setting $j=p$ in the summand in (\ref{Laguerre2}),
\begin{eqnarray}
\frac{(2p)!}{2} \cdot \frac{(-1)^{p+p}}{(2p)!} {2p \choose p} \varphi^{(k+p)}(0)\varphi^{(2p+k-p)}(0)
& = & \frac{1}{2}{2p \choose p} \gamma^2_{k+p} \\
& = & {2p-1 \choose p} \gamma_{k+p}^2\; .
\end{eqnarray}
For a fixed $j = 1, 2, \ldots, p$, and an arbitrary positive integer $p$, the coefficient of $\gamma_{k+p-j}\gamma_{k+p+j}$ is obtained by setting $j = p-j$ or $j = p+j$ in the summand in (\ref{Laguerre2}).  Thus, using the symmetry ${2p \choose p-j} = {2p \choose p+j}$,
\begin{equation}\hspace{-2.75 cm}
2\cdot \frac{(2p)!}{2}\cdot\frac{(-1)^{2p+j}}{(2p)!}{2p \choose p-j} \varphi^{(p+k-j)}(0)\varphi^{(p+k+j)}(0) 
\end{equation}
\begin{eqnarray}\hspace{7.15 cm}
& = & (-1)^j{2p \choose p-j} \gamma_{p+k-j}\gamma_{p+k+j} \; .
	\end{eqnarray}
\end{Proof}

\begin{Example}\label{L_k} 
Let $\varphi(x) = \sum_{k=0}^\infty \frac{\gamma_k}{k!} x^k$ be a real entire function.  The first 5 of the \textit{extended Tur\'an expressions} (cf. Theorem \ref{TLaguerre}) are: \\
\begin{tabular}{rl}
$\left. L_0(\varphi(x))\right|_{x=0}$ &= \; $\gamma_0^2$ \;, \\
$\left. L_1(\varphi(x))\right|_{x=0}$ &= \; $\gamma_1^2 - \gamma_0\gamma_2$ \;,\\
$\left. 12L_2(\varphi(x))\right|_{x=0}$ &= \; $3\gamma_2^2 - 4\gamma_1\gamma_3 + \gamma_0\gamma_4$ \;,\\
$\left. 360L_3(\varphi(x))\right|_{x=0}$ &= \; $10\gamma_3^2 - 15\gamma_2\gamma_4 + 6\gamma_1\gamma_5 - \gamma_0\gamma_6$ \;,  \\
$\left. 20160L_4(\varphi(x))\right|_{x=0}$ &= \; $35\gamma_4^2 - 56\gamma_3\gamma_5 + 28\gamma_2\gamma_6  -8\gamma_1\gamma_7 + \gamma_0\gamma_8$ \;.
\end{tabular}
\end{Example}

With the above coefficients, we define a class of non-linear operators that extend the non-linear operator of Theorem \ref{Stanley}.

\begin{Definition}\label{L^k} 
Let $\psi(x) = \sum_{k=0}^\infty a_k x^k$ be a real entire function.  For non-negative integers $p$, define the non-linear operators $a_k \mapsto L_k^p$\! , where $L_k^0 := a_k^2$\,, and for $p = 1, 2, 3, \ldots$ , set
\begin{equation}
L_k^p := {2p-1 \choose p} a_k^2 + \sum_{j=1}^p (-1)^j {2p \choose p-j} a_{k-j}a_{k+j} \; . \end{equation}
\end{Definition}

\vspace{0.5 cm}
\begin{Example} Let $\psi(x) = \sum_{k=0}^\infty a_k x^k$ be a real entire function.  The first 5 of the non-linear operators in Definition \ref{L^k} are: \\
$L_k^0 : a_k \mapsto a_k^2$ \;,\\
$L_k^1 : a_k \mapsto a_k^2 - a_{k-1}a_{k+1}$ \;,\\
$L_k^2 : a_k \mapsto 3a_k^2 - 4a_{k-1}a_{k+1} + a_{k-2}a_{k+2}$ \;,\\
$L_k^3 : a_k \mapsto 10a_k^2 - 15a_{k-1}a_{k+1} + 6a_{k-2}a_{k+2} - a_{k-3}a_{k+3}$ \;,\\
$L_k^4 : a_k \mapsto 35a_k^2 - 56a_{k-1}a_{k+1} + 28a_{k-2}a_{k+2} -8a_{k-3}a_{k+3} + a_{k-4}a_{k+4}$ \;.

\end{Example}

\begin{Notation}
In the sequel, we will allow $L_k^p$ to denote the non-linear operator $a_k \mapsto L_k^p$, and we will write $L_k^p\left[\psi(x)\right] = \sum_{k=0}^{n}L_k^p x^k$ to indicate the action of the non-linear operator $L_k^p$ on $\psi(x)= \sum_{k=0}^n a_kx^k$.
\end{Notation}

\subsection*{Symmetric Function Identities.}
For a fixed positive integer $n$ and $k = 1, 2, \ldots, n$, the elementary symmetric functions in the variables $z_1, z_2, \ldots, z_n$\,, are:
\begin{equation} \label{e_k}
e_k(z_1, z_2, \ldots, z_n) = \sum_{1\leq m_1 \leq m_2 \leq \cdots \leq m_k \leq k} \left( \prod_{j=1}^k z_{m_j} \right), 
\end{equation}
where $e_k(z_1, z_2, \ldots, z_n) = 0$, whenever $k \not\in \{1, 2, \ldots, n\},$ and $e_0(z_1, z_2, \ldots, z_n):=1$.

\vspace{0.5cm}
The proof of Theorem \ref{Stanley} depends on the properties of the non-linear operator $T_{\mu} : \Com[z] \to \Com[z]$,
\begin{equation}\label{T-op}
\sum_{k=0}^n a_k z^k \mapsto \sum_{i\leq j} \mu_{j-i}a_ia_jz^{i+j} \; ,
\end{equation}
and the following result.

\begin{Theorem}[\cite{BR10}]\label{Branden}
Let $\mu = \{\mu_k\}_{k=0}^{\infty}$ be a sequence of complex numbers and for \\ $k = 0, 1, \ldots, n$, let $e_k(z_1, z_2, \ldots, z_n)$ be the elementary symmetric functions in the variables $z_1, z_2, \ldots, z_n$.
Then,
\begin{equation}
\hspace{-3.5 cm} \sum_{i\leq j} \mu_{j-i} e_i(z_1, z_2, \ldots, z_n)e_j(z_1, z_2, \ldots, z_n) \end{equation} 
\[\hspace{2.75 cm} = e_n(z_1, z_2, \ldots, z_n)\sum_{k=0}^n \gamma_k e_{n-k}\left(z_1 + \frac{1}{z_1},z_2 + \frac{1}{z_2},  \ldots, z_n + \frac{1}{z_n}\right) \;, \]
where 
\begin{equation}\label{gamma_k} 
\gamma_k := \sum_{j=0}^{\lfloor \frac{k}{2} \rfloor} {k \choose j} \mu_{k-2j} \; .
\end{equation}
\end{Theorem}

\begin{Example}[\cite{BR10}]\label{mu app}
An application of Theorem \ref{Branden} with $\mu = \{1, 0, -1, 0, 0, \ldots \}$, the sequence obtained from the coefficients of $L_k^1$, yields
\begin{equation}\label{case1}
\sum_{k=0}^n \left( e^2_k(z_1, z_2, \ldots, z_n) - e_{k-1}(z_1, z_2, \ldots, z_n)e_{k+1}(z_1, z_2, \ldots, z_n) \right) 
\end{equation}
\[\hspace{2.5 cm} = e_n(z_1, z_2, \ldots, z_n) \sum_{k=0}^{\lfloor \frac{n}{2} \rfloor} C_k e_{n-2k} \left(z_1 + \frac{1}{z_1},z_2 + \frac{1}{z_2},  \ldots, z_n + \frac{1}{z_n}\right) \;,  \]
where $C_k := \frac{{2k\choose k}}{k+1}$ is the $k^{th}$ \textit{Catalan number}.
\end{Example}

We extend the identity in equation (\ref{case1}) to sequences $\mu$ obtained from the coefficients of an arbitrary $L_k^p$ (cf. Definition \ref{L^k}).

\begin{Lemma}\label{Q^p}
For fixed positive integers $n$ and $p$ define 
\begin{equation}
\hspace{-2.75 cm}  \mathcal{L}_k^p(z_1, z_2, \ldots, z_n) := {2p-1 \choose p}e^2_k(z_1, z_2, \ldots, z_n) \end{equation}
\[ \hspace{3.75 cm} + \sum_{j=1}^p (-1)^j{2p \choose p-j}e_{k+ j}(z_1, z_2, \ldots, z_n)e_{k- j}(z_1, z_2, \ldots, z_n) \] 
and 
\begin{equation}S(p,k) := \frac{{2p\choose p}{2k \choose k}}{{p+k\choose p}}\; .\end{equation}
Then, for a fixed positive integer $p$ and for any positive integer $n$,
\begin{equation} 
\hspace{-6.75 cm}\sum_{k=0}^n \mathcal{L}_k^p(z_1, z_2, \ldots, z_n) \end{equation}
\[ \hspace {1.75 cm} = e_n(z_1, z_2, \ldots, z_n)\sum_{k=0}^{\left\lfloor \frac{n}{2} \right\rfloor} \frac{S(p,k)}{2} e_{n-2k}\left(z_1 + \frac{1}{z_1},z_2 + \frac{1}{z_2},  \ldots, z_n + \frac{1}{z_n}\right) \;.\]  
\end{Lemma}

\begin{Proof}
Fix a positive integer $p$ and let the coefficients of $L_k^p$ (cf. Definition \ref{L^k}) form the the sequence $\{\mu_j\}_{j=0}^p$ by setting: $\mu_0 := {2p-1 \choose p}, \ \mu_{2j} := (-1)^j{2p \choose p-j}, \ \mu_{2j-1} := 0$, where $j = 1, 2, \ldots, p$.  Fix a positive integer $n$, and for integers $k$, where $0 \leq k \leq n$, define
\begin{equation} 
\gamma_k := \gamma_k(n,p) = \sum_{j=0}^{\lfloor \frac{k}{2} \rfloor} {k \choose j} \mu_{k-2j} \;.
\end{equation}
If $k$ is an odd integer, then $\gamma_k = 0$ and 
\begin{equation}
\label{gamma2}
\gamma_{2k} = \sum_{j=0}^k {2k\choose j} \mu_{2k-2j} = {2k\choose k} \mu_0 + \sum_{j=0}^{k-1} {2k \choose j} \mu_{2k-2j} \; .
\end{equation} 
By Theorem \ref{Branden}, it suffices to show that $\gamma_{2k} = \frac{S(p,k)}{2}$.  Using the symmetry of ${2k\choose j}$, we reverse the order of summation and obtain
\begin{eqnarray}
2\gamma_{2k} & = & 2{2k\choose k} \mu_0 + 2\sum_{j=0}^{k-1} {2k \choose j} \mu_{2k-2j} \\
             & = & 2{2k\choose k} \mu_0 + 2\sum_{j=1}^{k} {2k \choose k-j} \mu_{2j} \\
             & = & \sum_{j=-k}^{k} (-1)^j{2k \choose k-j} {2p\choose p-j} \; .
          \end{eqnarray}
\newpage
\noindent          
To complete the proof, we recall the following formula of Szily (\cite{Car}, \cite{Szily})
\begin{equation}
\label{szily}
\sum_{r = -b}^b (-1)^r {2a\choose a-r} {2b \choose b-r} = \frac{{2a \choose a}{2b \choose b}}{{a+b \choose a}} \; , 
\end{equation}
where $a \geq b$, and conclude that $2\gamma_{2k} = S(p,k)$.
\end{Proof}

\begin{Remarks} \noindent
\begin{enumerate}
\item $2C_k = S(1, k)$ and Example \ref{mu app} is equivalent to the assertion in Lemma \ref{Q^p} with $p = 1$.
\item $S(p,k)$ is an integer (cf. (\ref{szily}) and \cite{Cat}).
\item The numbers $S(p,k)$ are also known as the \textit{super Catalan numbers} (\cite{Cat}, \cite{Gessel}).

\end{enumerate}
\end{Remarks}

\subsection*{A Class of Hypergeometric Polynomials}  We now establish the properties of certain hypergeometric polynomials that appear in the proof of the main theorem.

\begin{Lemma} \label{hyper}
Let $p$ be a fixed positive integer.  Then, for any positive integer $n$, 
\begin{equation}
\label{hyper1}
\sum_{k=0}^{\left\lfloor \frac{n}{2} \right\rfloor} \frac{S(p,k)}{2} {n \choose 2k} z^k = {2p-1 \choose p} \ _2F_1\left( -\frac{n}{2}, \frac{1-n}{2}; p+1; 4z \right) \;, 
\end{equation}
where
\begin{equation}
S(p,k) = \frac{{2p\choose p}{2k \choose k}}{{p+k\choose p}}\; .
\end{equation}
\end{Lemma}

\begin{Proof}
Using the fact (\cite[Lemma 5, p. 22]{rain})
\begin{equation}
\label{fact1}
(\alpha)_{2k} = 2^{2k}\left( \frac{\alpha}{2} \right)_k \left( \frac{1+\alpha}{2} \right)_k \;,
\end{equation}
where $(\alpha)_k = \alpha(\alpha+1)(\alpha+2)\cdots (\alpha +k-1) = \frac{\Gamma(\alpha +k)}{\Gamma(\alpha)}$ is the Pochhammer symbol or ascending factorial, we rewrite the right member of (\ref{hyper1}):
\begin{eqnarray}
\hspace{1 cm} {2p-1\choose p} \sum_{k=0}^{\infty} \frac{\left( -\frac{n}{2} \right)_k \left( \frac{1-n}{2} \right)_k }{k!(p+1)_k} (4z)^k 
& = & \frac{1}{2}{2p\choose p} \sum_{k=0}^{\infty} \frac{(-n)_{2k}}{2^{2k}k!(p+1)_k} (4z)^k \\
                                                                              & = & \frac{1}{2}{2p\choose p} \sum_{k=0}^{\left\lfloor \frac{n}{2} \right\rfloor} \frac{(-n)_{2k}}{k!}\cdot \frac{1}{(p+1)_k}z^k \\
                                                                              & = & \frac{1}{2}{2p\choose p} \sum_{k=0}^{\left\lfloor \frac{n}{2} \right\rfloor} \frac{n!}{k!(n-2k)!} \cdot \frac{p!}{(p+k)!}z^k \\
                                                                              & = & \frac{1}{2}{2p\choose p} \sum_{k=0}^{\left\lfloor \frac{n}{2} \right\rfloor} {2k\choose k}{n\choose 2k} \cdot \frac{p!k!}{(p+k)!} z^k\\																											& = & \sum_{k=0}^{\left\lfloor \frac{n}{2} \right\rfloor} \frac{S(p,k)}{2} {n \choose 2k} z^k \; .
                                              																\end{eqnarray}
                                              																\end{Proof}

\begin{Lemma}\label{P^p}
For fixed positive integers $n$ and $p$, the zeros of the polynomial 
\begin{equation}
\label{poly}
Q^p_n(z) := \sum_{k=0}^{\left\lfloor \frac{n}{2} \right\rfloor} \frac{S(p,k)}{2} {n \choose 2k} z^k \;, \mbox{ where } S(p,k) = \frac{{2p\choose p}{2k \choose k}}{{p+k\choose p}}\;, 
\end{equation}
are all real and negative. 
\end{Lemma}
\begin{Proof}
Applying Lemma \ref{hyper}, we obtain
\begin{equation}
\label{original}
 Q^p_n(z) = {2p-1 \choose p} \ _2F_1\left( -\frac{n}{2}, \frac{1-n}{2}; p+1; 4z \right) \;.
\end{equation} 
We recall a formula relating the hypergeometric function $_2F_1$ and the Jacobi polynomial $P_n^{(\alpha, \beta)}(x)$ (\cite[formula (2),  p.254]{rain}):

\begin{equation}
\label{rainville}P_n^{(\alpha, \beta)}(x) = \frac{(1+\alpha)_n}{n!} \cdot \left( \frac{1+x}{2}\right)^n\! _2F_1\left( -n, -\beta -n; \alpha + 1; \frac{x-1}{x+1}\right) \;.
\end{equation}
If $n$ is an even integer, we let $n = 2m, \ m = 1, 2, \ldots$ , so that the right member of (\ref{original}) becomes

\begin{equation}
\label{neven}
{2p-1 \choose p} \ _2F_1\left(-m, \frac{1-2m}{2}; p+1; 4z \right) \; .
\end{equation}
Setting $\alpha := p, n := m, \beta := -\frac{1}{2}, x := \frac{1+4z}{1-4z}$ in (\ref{rainville}), yields

\begin{equation}\label{Jacobi1} 
P_{m}^{(p, -\frac{1}{2})}\left( \frac{1+4z}{1-4z} \right) = \frac{(1+p)_{m}}{m!} \cdot \left( z-\frac{1}{4}\right)^{m}\! _2F_1\left( -m, \frac{1-2m}{2}; p+1; 4z\right) \;. 
\end{equation}
If $n$ is an odd integer, we let $n = 2m-1, m=1,2, \ldots$ , so that the right member of (\ref{original}) becomes

\begin{equation}
\label{nodd}
{2p-1 \choose p} \ _2F_1\left(1-m, \frac{1-2m}{2}; p+1; 4z \right) \; .
\end{equation}
Setting $\alpha := p, n := m-1, \beta := \frac{1}{2}, x := \frac{1+4z}{1-4z}$ in (\ref{rainville}), we obtain

\begin{equation}\label{Jacobi2}
P_{m-1}^{(p, \frac{1}{2})}\left( \frac{1+4z}{1-4z} \right) = \frac{(1+p)_{m-1}}{(m-1)!} \cdot \left( z-\frac{1}{4}\right)^{m-1}\! _2F_1\left( 1-m, \frac{1-2m}{2}; p+1; 4z\right) \;. 
\end{equation}
It is well known that if $\alpha >-1$ and $\beta >-1$, then the zeros of the Jacobi polynomial $P_n^{(\alpha, \beta)}(x)$ are distinct and lie in the interval $(-1,1)$ (cf. \cite[p. 261]{rain}).  Thus, if $\gamma \in (-1,1)$ is a zero of $P_m^{(p, -\frac{1}{2})}\left( \frac{1+4z}{1-4z} \right)$ or a zero of $P_{m-1}^{(p, \frac{1}{2})}\left( \frac{1+4z}{1-4z} \right)$, then a calculation shows that $z = \frac{\gamma -1}{4(\gamma+1)} < 0$.  Therefore, the zeros of the polynomials in equations (\ref{Jacobi1}) and (\ref{Jacobi2}), and whence the zeros of the polynomial $Q_n^p(z)$, are all real and negative.
\end{Proof}

\subsection*{Supplementary Results}  
\begin{Theorem}[Grace-Walsh-Szeg\"o, \cite{ATP}]
\label{GWS}
Let $f \in \Com[z_1, z_2, \ldots, z_n]$ be a multi-affine and symmetric polynomial.  Let $K$ be a circular region and suppose that either $K$ is convex or that the degree of $f$ is $n$.  Then, for any $\zeta_1, \zeta_2, \ldots, \zeta_n \in K$, there exists a $\zeta \in K$ such that $f(\zeta_1, \zeta_2, \ldots, \zeta_n) = f(\zeta, \zeta, \ldots, \zeta)$.
\end{Theorem}

\begin{Theorem}[Sz\'asz, \cite{Szasz}] \label{Szasz} Let $\{P_n(z)\}_{n=0}^{\infty}$ be a sequence of polynomials, where \\ $P_n(z) = \sum_{k=0}^{M_n} a_{n,k}z^k, a_{n,0} \neq 0, a_{n,M_n} \neq 0, M_n \to \infty$, all of whose zeros lie in an open half-plane $H \subset \Com$ with boundary containing the origin.  If for some constants $\alpha_0, \alpha_1,$ and all $n$,
\begin{equation}\label{szasz bdds}
0 < \alpha_0 \leq |a_{n,0}| \leq \alpha_1, \ |a_{n,1}| \leq \alpha_1, \ |a_{n,2}| \leq \alpha_1 < \infty \; ,
\end{equation}
then the sequence $\{P_n(z)\}_{n=0}^{\infty}$ is uniformly bounded in any circle $|z| \leq r$, and in fact, 
\begin{equation}
|P_n(z)| \leq \alpha_1 \exp \left( r \frac{\alpha_1}{\alpha_2} + 3r^2 \left( \frac{\alpha_1^2}{\alpha_0^2} + \frac{\alpha_1}{\alpha_0}\right) \right) \;.
\end{equation}
\end{Theorem}
\section{Main Theorem}
The analog of Theorem \ref{Stanley} (cf. \cite[proof of Conjecture 1.1]{BR10}) for the extended Tur\'an inequalities is that the operator $L_k^p$ (cf. Definition \ref{L^k}), for an arbitrary positive integer $p$, takes real polynomials with all real negative zeros into polynomials of the same type.

\vspace{0.5cm}
\begin{Theorem}\label{ExtStanley}Let $\psi(z) = \sum_{k=0}^n a_kz^k = \prod_{k=1}^n(1+\rho_k z)$, where $\rho_k >0$ for $1 \leq k \leq n$, be a real polynomial with all real negative zeros.  Let $p$ be a positive integer and let $L_k^p$ be the non-linear operator (cf. Definition \ref{L^k}) 
\begin{equation}\label{main3}
a_k \mapsto {2p-1 \choose p} a_k^2 + \sum_{j=1}^p (-1)^j {2p \choose p-j} a_{k-j}a_{k+j} \; .
\end{equation}
Then, the zeros of the polynomial
\begin{equation}
\label{main1}
L_k^p\left[\psi(z)\right] = \sum_{k=0}^n L_k^p z^k \end{equation}
are all real and negative.
\end{Theorem}

\begin{Proof}
Fix a positive integer $p$ and suppose that there exists a $\xi$ with $\R \xi >0$ such that $\xi^2 = \zeta$ and $L_k^p\left[\psi(\zeta)\right] = 0$.  From (\ref{main3}) and (\ref{main1}) we have

\begin{equation}   \hspace{-.25 cm} L_k^p\left[\psi(\zeta)\right] =  \sum_{k=0}^n \left( {2p-1 \choose p} a_k^2 + \sum_{j=1}^p (-1)^j {2p \choose p-j} a_{k-j}a_{k+j}\right)  \zeta^k
\end{equation}

\begin{equation}\label{main2}
\hspace{1.4 cm} = \sum_{k=0}^n \left( {2p-1 \choose p} a_k^2\xi^{2k} + \sum_{j=1}^p (-1)^j {2p \choose p-j} a_{k-j}a_{k+j}\xi^{2k}\right) \; .
\end{equation}
Using the properties of the elementary symmetric functions, we apply Lemma \ref{Q^p} to (\ref{main2}) and obtain

\begin{equation}
\hspace{-4.5cm} \sum_{k=0}^n \left( {2p-1 \choose p} e_k^2(\rho_1\xi, \rho_2\xi, \ldots, \rho_n\xi) \right.
\end{equation}

\[ \hspace{2 cm} + \left. \sum_{j=1}^p (-1)^j {2p \choose p-j} e_{k-j}(\rho_1\xi, \rho_2\xi, \ldots, \rho_n\xi)e_{k+j}(\rho_1\xi, \rho_2\xi, \ldots, \rho_n\xi) \right)
\]

\begin{equation}  = a_n\xi^n \sum_{k=0}^{\left\lfloor \frac{n}{2} \right\rfloor} \frac{S(p,k)}{2} e_{n-2k}\left(\rho_1\xi + \frac{1}{\rho_1\xi}, \ \rho_2\xi + \frac{1}{\rho_2\xi}, \  \ldots, \ \rho_n\xi + \frac{1}{\rho_n\xi}\right) \;.
\end{equation}
From the assumptions that $\R \xi >0$ and $\rho_k >0$ for $1 \leq k \leq n$, it follows that $\R (\rho_k\xi + \frac{1}{\rho_k\xi}) >0$ for $1\leq k \leq n$.  By the Grace-Walsh-Szeg\" o Theorem, there exists an $\eta \in \Com$ with $\R \eta >0$ such that
\begin{eqnarray}
0 & = & \sum_{k=0}^{\left\lfloor \frac{n}{2} \right\rfloor} \frac{S(p,k)}{2} e_{n-2k}\left(\eta, \eta, \ldots, \eta\right) \\
  & = & \sum_{k=0}^{\left\lfloor \frac{n}{2} \right\rfloor} \frac{S(p,k)}{2} {n \choose 2k} \eta^{n-2k} \\
  & = & \eta^n \sum_{k=0}^{\left\lfloor \frac{n}{2} \right\rfloor} \frac{S(p,k)}{2} {n \choose 2k} \left( \frac{1}{\eta^2}\right) ^k \\ 
  & = & \eta^n Q_n^p \left( \frac{1}{\eta^2} \right) \;,
  \end{eqnarray}
where $Q_n^p(z) = \sum_{k=0}^{\left\lfloor \frac{n}{2} \right\rfloor} \frac{S(p,k)}{2} {n \choose 2k} z^k$.  Since $\zeta = \xi^2$, we know that $\zeta, \frac{1}{\eta^2} \in \Com \sm \{x\in \Rea \; : \; x \leq 0\}$ and, therefore, the zeros of $\sum_{k=0}^n L_k^p z^k$ are all real and negative, provided that the zeros of $Q_n^p(z)$ are all real and negative for a fixed positive integer $p$ and any positive integer $n$.  This is the assertion in Lemma \ref{P^p}.
\end{Proof}

\vspace{0.5 cm}
In a sequel to the proof of Theorem \ref{Branden}, Br\"and\'en gives several equivalent stability results concerning the action of the operator $T_{\mu}$ (cf. (\ref{T-op})) on polynomials and transcendental entire functions (cf. \cite[Theorems 3.2 and 5.4]{BR10}):
\begin{enumerate}
\item $T_{\mu} \left((1 +z)^n \right) \neq 0$, whenever $\R z >0$; 
\item For all polynomials $p(z)$ with only real and negative zeros, the polynomial $T_{\mu}(p(z))$ is either identically zero or $T_{\mu}(p(z)) \neq 0$, whenever $\R z >0$; 
\item For all $\varphi(x) \in$ \LPP, the function $T_{\mu}(\varphi(x))$ is either identically zero or the uniform limit, on compact subsets of $\Com$, of polynomials all of whose zeros lie in the closed left half-plane.
\end{enumerate}

The class of non-linear operators $L_k^p$ enjoys the above properties, because the coefficients of an arbitrary $L_k^p$ define the sequence $\mu$ as in Lemma \ref{Q^p}, and by Theorem \ref{ExtStanley} satisfy statements (i) and (ii) above.  We refine statement (iii) and show that for this choice of $\mu$ the non-linear operators $L_k^p$ preserve \LPP.

\subsection*{A Transcendental Characterization}
\begin{Theorem}\label{trans} For a fixed positive integer $p$, let $L_k^p$ be the non-linear operator
\begin{equation}
a_k \mapsto {2p-1 \choose p} a_k^2 + \sum_{j=1}^p (-1)^j {2p \choose p-j} a_{k-j}a_{k+j} \; ,
\end{equation}
and let $\varphi(z) = \sum_{k=0}^{\infty} \frac{\gamma_k}{k!}z^k  = \sum_{k=0}^{\infty} a_k z^k \in$ \LPP \!. Then, $L_k^p[ \varphi(z) ] \in$ \LPP .
\end{Theorem}

\begin{Proof}
Fix a positive integer $p$.  To establish $L_k^p[ \varphi(z) ] \in$ \LPP, it suffices to approximate $L_k^p[ \varphi(z) ]$, uniformly on compact subsets of $\Com$, by real polynomials all of whose zeros are real and of the same sign (cf. \cite[Lemma 2.2]{CC89}).  The $n^{th}$ \textit{Jensen polynomial} associated with $\varphi(z)$, given by $g_n(z) := \sum_{k=0}^n {n\choose k} \gamma_k z^k, n = 0, 1, \ldots$ , has only real negative zeros (cf. \cite{PS}).  The zeros of $g_n(\frac{z}{n}), n = 1, 2, \ldots$\, ,  are all real and negative and by Theorem \ref{ExtStanley}, for a fixed positive integer $p$, the zeros of the polynomial $L_k^p[ g_n(\frac{z}{n}) ]$ are all real and negative.  For $n \geq 2$,
\begin{equation}
g_n \left(\frac{z}{n} \right) = \gamma_0 + \gamma_1 z + \sum_{k=2}^n {n \choose k} \gamma_k \left( \frac{z}{n} \right)^k 
\end{equation}
and
\begin{equation}
L_k^p \left[ g_n\left(\frac{z}{n} \right) \right] = L_0^p + L_1^pz + \sum_{k=2}^n L_k^p z^k \; .
\end{equation}
For each $k \in \Nat$,
\begin{equation} \label{lim1}
\lim_{n\to \infty} {n \choose k} \frac{\gamma_k}{n^k} = 
\frac{\gamma_k}{k!} \lim_{n\to \infty} \prod_{r=1}^{k-1} \left( 1-\frac{r}{n} \right) = \frac{\gamma_k}{k!} \;,
\end{equation}
and a calculation using (\ref{lim1}), yields
\begin{equation} \hspace{0 cm}
\lim_{n\to \infty} \ L_k^p\left[ \sum_{k=2}^n {n \choose k} \gamma_k \left( \frac{z}{n} \right)^k \right] 
= L_k^p \left[ \sum_{k=2}^{\infty} \frac{\gamma_k}{k!} z^k \right]\;.
\end{equation}

\noindent
Let $\{P_n(z)\}_{n=2}^{\infty}$ be the family of polynomials $P_n(z) := L_k^p \left[ g_n\left(\frac{z}{n} \right) \right] = \sum_{k=0}^n a_{n,k}z^k$.  By Theorem \ref{ExtStanley},
\begin{equation}
\begin{tabular}{lcl} 
$a_{n,0}$ &=& ${2p-1 \choose p} \gamma_0^2$ \;, \\
$a_{n,1}$ &=& ${2p-1 \choose p} \gamma_1^2 - {2p \choose p-1}{n \choose 2}\frac{\gamma_0\gamma_2}{n^2}$ \;, \\
$a_{n,2}$ &=& ${2p-1 \choose p} \left( {n\choose 2}\frac{\gamma_2}{n^2} \right)^2 - {2p \choose p-1}{n \choose 3}\frac{\gamma_1\gamma_3}{n^3} + {2p \choose p-2}{n \choose 4}\frac{\gamma_0\gamma_4}{n^4}$  \;,
\end{tabular}
\end{equation}
are non-negative for any positive integer $p$.  Therefore, there exist positive constants $\alpha_0, \alpha_1$, satisfying (\ref{szasz bdds}).  By Theorem \ref{Szasz}, $\{P_n(z)\}_{n=2}^{\infty}$ is locally uniformly bounded on compact subsets of $\Com$, and by well-known results (see for example \cite[p.333]{levin} or Montel's Theorem \cite[pp.\,21-30]{Montel}), there exists a subsequence of $\{P_n(z)\}_{n=2}^{\infty}$ converging, uniformly on compact subsets of $\Com$, to the entire function $L_k^p \left[  \varphi(z) \right]$.  

\end{Proof}

\vspace{0.5 cm}
\section{Applications and Examples}
\begin{Example}
The analog of Theorem \ref{Stanley} fails for the \textit{extended Tur\'an expressions} (cf. Example \ref{L_k}).  Consider the polynomial $\psi(x) = (1+x)^3 = 1 + 3x + 3x^2 + x^3 = \sum_{k=0}^3 \frac{\gamma_k}{k!} x^k$.  The polynomial 
\begin{equation} \sum_{k=0}^3(3\gamma_k^2 - 4\gamma_{k-1}\gamma_{k+1} + \gamma_{k-2}\gamma_{k+2})x^k = 12 + 84x + 36x^2 + 108x^3
\end{equation}
has 2 non-real zeros.
\end{Example}

\begin{Example}
The hypothesis on the degrees of both the polynomials in Theorem \ref{ExtStanley} (cf. Theorem \ref{Stanley}) cannot be relaxed.  Consider the polynomial $\psi(x) = (1+x)^5 =  1 + 5x + 10x^2 + 10x^3 + 5x^4 + x^5 = \sum_{k=0}^5 a_k x^k$.  The polynomial
\begin{equation}
L_k^2\left[\psi(x)\right] = \sum_{k=0}^5 (3a_k^2 - 4a_{k-1}a_{k+1} + a_{k-2}a_{k+2}) x^k = 3 + 35x + 105x^2 + 105x^3 + 35x^4 + 3x^5
\end{equation}
has all negative real zeros, but the polynomial
\begin{equation}
\sum_{k=0}^4 (3a_k^2 - 4a_{k-1}a_{k+1} + a_{k-2}a_{k+2}) x^k = 3 + 35x + 105x^2 + 105x^3 + 35x^4
\end{equation}
has a pair of non-real zeros.  
\end{Example}

\subsection*{A Question of Fisk}Fisk asked whether the class of non-linear operators $a_k \mapsto S_r := a_k^2 -  a_{k-r}a_{k+r}, r = 0, 1, 2, \ldots$ , acting on functions of the form $\sum_{k=0}^{n} a_kx^k$, takes polynomials with only real negative zeros into polynomials of the same type (\cite[Question 2]{Fisk}).  Here, as with $L_k^p$, we set $a_k = 0$ for $k<0$ (cf. Notation \ref{not1}).

$S_0$ produces the zero polynomial and Br\"and\'en (\cite{BR10}) confirmed the cases when $r = 1, 2, 3$.  Recently, R. Yoshida (\cite{Y10}) confirmed the case when $r = 4$ and produced a counterexample in the case $r = 6$.  Nevertheless, the non-linear operators $S_r$ and $L_k^p$ are related in a remarkable way.

\begin{Proposition}\label{LS}
Fix a positive integer $p$ and let $L_k^p$ be the non-linear operator of Definition \ref{L^k}.  
Then, 
\begin{equation}\label{S^r}
L_k^p = \sum_{j=1}^p (-1)^{j+1}{2p \choose p-j}(a_k^2 -  a_{k-j}a_{k+j}) \; .
\end{equation}
\end{Proposition}

\begin{Proof}Using (\ref{szily}) with $a = 0$, we rewrite the right member of (\ref{S^r}) and obtain
\begin{eqnarray}
L_k^p & = & a_k^2 \sum_{j=1}^p (-1)^{j+1}{2p \choose p-j} + \sum_{j=1}^p (-1)^j{2p \choose p-j} a_{k-j}a_{k+j}\\
      & = & {2p-1 \choose p} a_k^2 + \sum_{j=1}^p (-1)^j {2p \choose p-j} a_{k-j}a_{k+j} \; .
      \end{eqnarray}     
\end{Proof}

\begin{Proposition}\label{FiskCon} If the zeros of the real polynomial $\varphi(x) = \sum_{k=0}^n a_kx^k$ are all real and negative, then for any positive integer $p$, the zeros of the polynomial
\begin{equation}\label{Lk and Sk}
\sum_{j=1}^p \left( (-1)^{j+1}{2p \choose p-j} \sum_{k=0}^n (a_k^2 -  a_{k-j}a_{k+j}) x^k \right)
\end{equation}
are all real and negative.
\end{Proposition}

\begin{Proof}
Suppose that the zeros of the polynomial $\varphi(x) = \sum_{k=0}^n a_kx^k$ are all real and negative.  By Theorem \ref{ExtStanley}, the zeros of the polynomial $L_k^p\left[\varphi(x)\right]$ are all real and negative for any positive integer $p$, and by Proposition \ref{LS}
\begin{eqnarray}
L_k^p\left[\varphi(x)\right] & = & \sum_{k=0}^n \left( \sum_{j=1}^p (-1)^{j+1}{2p \choose p-j} (a_k^2 -  a_{k-j}a_{k+j}) \right) x^k \\
                             & = & \sum_{j=1}^p \left( (-1)^{j+1}{2p \choose p-j} \sum_{k=0}^n (a_k^2 -  a_{k-j}a_{k+j}) x^k \right) \; .
               							 \end{eqnarray}
\end{Proof}

\subsection*{Totally Positive Sequences} According to the theory of totally positive sequences developed by Aissen, Edrei, Schoenberg, and Whitney (\cite{ASW}, \cite{Edrei}), the coefficients $\{a_k\}_{k=0}^{\infty}$ of a function $\varphi(x) = \sum_{k=0}^{\infty} a_kx^k \in$ \LPP form a totally positive sequence.  Therefore, by Theorem \ref{ExtStanley} the sequences $\{L_k^p\}_{k=0}^{\infty}$ are totally positive sequences for each positive integer $p$, provided that the non-linear operators defined by $L_k^p$ are acting on functions in \LPP.

For positive integers $r$, the non-linear operator $S_r$, in the determinant form 
\begin{equation}
a_k \mapsto \left| \begin{tabular}{cc} $a_k$ & $a_{k-r}$ \\ $a_{k+r}$ & $a_k$ \end{tabular} \right| \;,
\end{equation} 
occurs naturally as a minor of the 4-way infinite Toeplitz matrix $A = (a_{ij})$, obtained from the coefficients of $\varphi(x) = \sum_{k=0}^{\infty} a_kx^k \in$ \LPP, by setting $a_{ij} = a_{k+(i-j)}$:
\begin{equation}\label{toeplitz} A= 
\left(
\begin{tabular}{ccccccc}
         & $\vdots$ & $\vdots$  & $\vdots$  &          & $\vdots$  &  \\
$\cdots$ & $a_k$    & $a_{k-1}$ & $a_{k-2}$ & $\cdots$ & $a_{k-r}$ & $\cdots$ \\
$\cdots$ & $a_{k+1}$    & $a_{k}$ & $a_{k-1}$ & $\cdots$ & $a_{k-r+1}$ & $\cdots$ \\
$\cdots$ & $a_{k+2}$    & $a_{k+1}$ & $a_{k}$ & $\cdots$ & $a_{k-r+2}$ & $\cdots$ \\
         & $\vdots$ & $\vdots$  & $\vdots$  & $\ddots$ & $\vdots$  &  \\
$\cdots$ & $a_{k+r}$    & $a_{k+r-1}$ & $a_{k+r-2}$ & $\cdots$ & $a_{k}$ & $\cdots$ \\
         & $\vdots$ & $\vdots$  & $\vdots$  &          & $\vdots$  & 
\end{tabular}
\right) \; .
\end{equation}

With the exception of $L_k^1 = S_1$, it seems that the non-linear operators defined by $L_k^p$ cannot be realized as a minor of the above matrix $A$.  From (\ref{S^r}) and the Laplace expansion along the $i^{th}$ row of the determinant of a $p\times p$ submatrix $B:= (b_{ij})$ of $A$, 
\begin{equation}
\det B = \sum_{j=1}^p (-1)^{j+i} b_{ij} \det B(i|j) \; ,
\end{equation}
where $B(i|j)$ is the submatrix of $B$ obtained by deleting row $i$ and column $j$, it is clear that $L_k^p$ is not the determinant of any submatrix $B$ of $A$ larger than $2 \times 2$.  Moreover, for each positive integer $p$, the non-linear operator $L_k^p$ consists of the term ${2p-1 \choose p} a_k^2$ and $p$ other terms, and for $p >1$ cannot be realized as a $2 \times 2$ determinant, and \textit{a fortiori}, any minor of $A$ larger than $2 \times 2$.

\section*{Acknowledgments}
The author is indebted to Dr. George Csordas, Mr. Matthew Chasse, and Mr. Rintaro Yoshida for their careful reading of the manuscript and many helpful suggestions.





\end{document}